\documentclass[10pt]{amsart}
\usepackage{amstext,amsfonts,amssymb,amscd,amsbsy,amsmath}
\usepackage{ifthen}
\usepackage{amsthm}
\usepackage{latexsym}
\usepackage[all]{xy} \SelectTips{eu}{} \SilentMatrices 
\usepackage{enumerate}
\usepackage{hyperref}

\def\Pic{\operatorname{Pic}}

\newcommand{\tensor}{\otimes}
\newcommand{\iso}{\cong}
\newcommand{\onto}{\to\hspace*{-.8em}\to}
\newcommand{\into}{\hookrightarrow}

\newcommand{\dm}{\operatorname{dim}}

\renewcommand{\dim}[2][]{\operatorname{dim}_{#1}(#2)}

\newcommand{\G}[1]{G(#1)}
\newcommand{\K}[1]{K(#1)}
\renewcommand{\H}[4][]{H_{#1}^{#2}(#3, #4)}
\newcommand{\Hp}[2][]{\operatorname{P}_{#1}(#2)}
\newcommand{\Hpi}[3][]{\operatorname{P}_{#1}^{(#2)}(#3)}
\newcommand{\Hs}[2][]{\operatorname{H}_{#1}(#2)}

\newcommand{\ip}[1]{\langle #1 \rangle}
\newcommand{\length}{\operatorname{length}}
\newcommand{\Proj}[1]{{\operatorname{Proj}}(#1)}
\newcommand{\Spec}[1]{{\operatorname{Spec}}(#1)}
\newcommand{\Tor}{\operatorname{Tor}}

\newcommand{\ds}{\displaystyle}

\newcommand{\sint}{\mbox{$\int$}}

\newcommand{\bC}{\mathbb{C}}
\newcommand{\C}{\mathbb{C}}

\newcommand{\cF}{\mathcal{F}}

\newcommand{\sk}{k}

\newcommand{\bN}{\mathbb{N}}
\newcommand{\fm}{\mathfrak{m}}

\newcommand{\cO}{{\mathcal{O}}}
\newcommand{\fp}{\mathfrak{p}}
\newcommand{\bQ}{\mathbb{Q}}
\newcommand{\Q}{\mathbb{Q}}
\newcommand{\bQl}{\mathbb{Q}_\ell}
\newcommand{\bZl}{\mathbb{Z}_\ell}

\newcommand{\cT}{\mathcal{T}}
\newcommand{\bZ}{\mathbb{Z}}
\newcommand{\Z}{\mathbb{Z}}

\newcommand{\im}{\operatorname{im}}
\newcommand{\ch}[1][]{ch_{#1}}
\newcommand{\chj}[2][]{\operatorname{ch}_{#1}^{#2}}
\newcommand{\CH}[2][]{CH^{#1}(#2)}
\newcommand{\Td}[1]{\operatorname{Td}(#1)}

\newcommand{\cy}[1][]{{cy}_{#1}}

\newcommand{\conjref}[1]{Conjecture \ref{#1}}
\newcommand{\corref}[1]{Corollary \ref{#1}}
\newcommand{\eqnref}[1]{(\ref{#1})}

\newcommand{\lemref}[1]{Lemma \ref{#1}}

\newcommand{\thmref}[1]{Theorem \ref{#1}}

\def\lra{\longrightarrow}

\def\map#1{{\buildrel #1 \over \lra}} 
\def\bu{\bullet}

\def\cO{\mathcal O}

\def\cF{\mathcal F}

\newcommand{\bP}{\mathbb{P}}

\def\chr{\operatorname{char}}

\def\lra{\longrightarrow}

\newcommand{\et}{{\operatorname{\acute{e}t}}}
\newcommand{\topo}{{\operatorname{top}}}

\newcommand{\btt}{\begin{tt}}
\newcommand{\ett}{\end{tt}}

\newtheorem{theorem}{Theorem}[section]
\newtheorem{thm}[theorem]{Theorem}
\newtheorem{lem}[theorem]{Lemma}
\newtheorem{cor}[theorem]{Corollary} 
\newtheorem{conj}[theorem]{Conjecture}
\newtheorem{prop}[theorem]{Proposition}

\theoremstyle{definition}

\newtheorem{example}[theorem]{Example}

\newtheorem*{chunk*}{}

\numberwithin{equation}{section}

\theoremstyle{remark}
\newtheorem{rem}[theorem]{Remark}

\begin{document}


\title[Hochster's Theta invariant]
{Hochster's Theta invariant and the Hodge-Riemann bilinear relations}

\author{W. Frank Moore} \address{Department of Mathematics, Cornell
  University, Ithaca, NY 14853} \email{frankmoore@math.cornell.edu}

\author{Greg Piepmeyer} \address{Department of Mathematics, University
  of Missouri, Columbia, MO 65211} \email{piepmeyerg@missouri.edu}

\author{Sandra Spiroff} \address{Department of Mathematics, University
  of Mississippi, University, MS 38677} \email{spiroff@olemiss.edu}
\thanks{The third author was supported in part by AWM-NSF grant
  192332}

\author{Mark E. Walker} \address{Department of Mathematics, University
  of Nebraska, Lincoln, NE 68588} \email{mwalker5@math.unl.edu}
\thanks{The fourth author was supported in part by NSF grant DMS-0601666}

\begin{abstract}  
  Let $R$ be an isolated hypersurface singularity, and let $M$ and $N$
  be finitely generated $R$-modules.  As $R$ is a hypersurface, the
  torsion modules of $M$ against $N$ are eventually periodic of period
  two (i.e., $\Tor_i^R(M,N) \iso \Tor_{i+2}^R(M,N)$ for $i \gg
  0$). Since $R$ has only an isolated singularity, these torsion
  modules are of finite length for $i \gg 0$.  The theta invariant of
  the pair $(M,N)$ is defined by Hochster to be $\length
  (\Tor_{2i}^R(M,N)) - \length (\Tor_{2i+1}^R(M,N))$ for $i \gg 0$.
  H.\ Dao has conjectured that the theta invariant is zero for all
  pairs $(M,N)$ when $R$ has even dimension and contains a field.
  This paper proves this conjecture under the additional assumption
  that $R$ is graded with its irrelevant maximal ideal giving the
  isolated singularity.  We also give a careful analysis of the theta
  pairing when the dimension of $R$ is odd, and relate it to a
  classical pairing on the smooth variety $\Proj{R}$.
\end{abstract}

\maketitle

\section{introduction}

If $R$ is a hypersurface --- that is, a quotient of a regular ring $T$ by
a single element --- and $M$ and $N$ are finitely generated $R$-modules,
then the long exact sequence %
\begin{equation} \label{les}
  \cdots \! \to \! \Tor_n^T(M,N) \!\to\! \Tor_n^R(M,N) \!\to\!
  \Tor_{n-2}^R(M,N) \!\to\! \Tor_{n-1}^T(M,N) \! \to \! \cdots   
\end{equation}
coming from \cite[Chapter XV]{CE} shows that 
\begin{equation*}
\Tor_i^R(M,N) \iso \Tor_{i+2}^R(M,N)   \quad {\text{for $i \gg 0$.}}
\end{equation*}
When these torsion modules are of finite length, M.~Hochster \cite[Theorem 1.2]{Ho}
defines 
\begin{equation*}
  \theta^R(M,N) = 
  \length(\Tor^R_{2i}(M,N)) - \length(\Tor^R_{2i+1}(M,N))
  \quad {\text{for $i \gg 0$.}}
\end{equation*}
If $R$ has at most a finite number of singularities, then
$\theta^R(M,N)$ is defined for all pairs $(M,N)$ of finitely generated
$R$-modules.  Typically we write just $\theta$ for $\theta^R$.  

Hochster introduced the $\theta$ pairing in his
study of the Direct Summand Conjecture: if $A \subseteq B$ is a module
finite ring extension of a regular ring $A$, then $A$ is a direct
summand of the $A$-module $B$.  The Direct Summand Conjecture is known
if $A$ is equicharacteristic \cite[Theorem 2]{Ho2} or has dimension at most three \cite{Heit}.  
Hochster showed that the Direct Summand Conjecture holds provided
$\theta(S/\fp, -)$ is the zero function for a particular prime $\fp
\in \Spec{S}$, where $S$ belongs to a explicit family of (mixed
characteristic) local hypersurfaces.  %

Hochster \cite[Theorem 1.4]{Ho} showed that if $R$ is an
{\emph{admissible}}%
\footnote{$R$ is {\emph{admissible}} if a completion $\hat{R}$ of $R$
  at a maximal ideal satisfies $\hat{R} \iso T/(f)$, and the
  {\emph{dimension inequality}}, {\emph{vanishing}}, and
  {\emph{positivity}} of Serre \cite[V.5.1]{Se} hold for $T$.  Serre
  showed that these conditions on $T$ hold when $T$ is a regular local
  ring containing a field.  } %
hypersurface and $\length(M \otimes_R N) < \infty$, then $\theta(M, N)
= 0$ if and only if $\dim{M} + \dim{N} \leq \dim{R}$.  H.~Dao
\cite{Dao1}, \cite{Dao2} studied the vanishing of $\theta$ for
admissible local hypersurfaces $R$ which have only isolated
singularities.  These papers motivated our work; in particular we
address the following conjecture of Dao:

\begin{conj} \cite[Conjecture 3.15]{Dao1} \label{Dao Conj} %
Let $R$ be an isolated  hypersurface singularity.  Assume that
$\dim{R}$ is even and $R$ contains a field.  Then $\theta(M,N) = 0$
for all pairs of finitely generated $R$-modules $M$ and $N$.
\end{conj}

In this paper, we prove this conjecture when $R$ is a graded, finitely
generated algebra over a field $\sk$ that is non-singular away from
its irrelevant maximal ideal; see \thmref{n even result} for our
precise statement.

The $\theta$ pairing induces a symmetric bilinear form on the
Grothendieck group of $R$.  When $n = \dim{R}$ is odd and $\sk$ is
separably closed, we prove $\theta$ factors through the Chern
character map taking values in \'etale cohomology; see \thmref{n odd thm}.
Moreover, when $\chr{k} = 0$, 
we show in \thmref{n odd pos defn}
that $(-1)^{\frac{n+1}{2}}\theta$ is positive semi-definite,
and when $\sk = \bC$, we identify its kernel using the Hodge-Riemann bilinear relations.

In section \ref{sec:nonstandardgrading} we extend our results on the
$\theta$ pairing (under the same assumptions on the ground field) to graded hypersurfaces
$S = k[y_1,\dots,y_n]/(g)$ where $\deg y_i = e_i \geq 1$, and $g$ is homogeneous
with respect to this grading.

Another source of interest in the $\theta$ pairing comes from a result
of Dao \cite[Proposition 2.8]{Dao1}, which provides a connection
between the vanishing of $\theta$ and the {\emph{rigidity}}%
\footnote{A pair of modules $(M, N)$ is {\emph{rigid}} if for any
  integer $i \geq 0$, $\Tor_i^R(M,N) = 0$ implies $\Tor_j^R(M,N) = 0$
  for all $j \geq i$.  A module $M$ is {\emph{rigid}} if for all $N$
  the pair $(M, N)$ is rigid.}  %
of $\Tor$.  Namely, when $R$ is an admissible hypersurface and $M, N$
are $R$-modules such that $\theta(M,N)$ is defined, then $\theta(M,N)=
0$ implies rigidity of the pair $(M,N)$.  Our results imply that, if
$R$ is a graded $\sk$-algebra as above, with $\chr{k} = 0$, then an $R$-module
$M$ is rigid if $\theta(M,M)=0$; see \corref{n odd theta rigidity}.  

Finally, for readers familiar with the Herbrand difference, perhaps
through \cite[Section 10.3]{Bu}, we note that the Herbrand difference
and $\theta$ are closely related. In detail, each can be interpreted
as coming from a pairing on the graded rational Chow group, $\CH[\bullet]{R}_\Q$,
of $R$. On the component $\CH[j]{R}_\Q$, they coincide for $j$ odd
and 
differ by a sign for $j$ even.  In particular, our results show that,
over a field of characteristic zero, the Herbrand difference is a negative semi-definite form;
see Example \ref{Buchweitz ex}.

\section{Background} \label{background}

Throughout the rest of this paper, we make the following assumptions:
\begin{equation} \label{assumptions}
  \begin{gathered}
    \begin{tabular}{l@{\hspace*{1em}}l}
      $\bu$ & $\sk$ is a field. \\
      \begin{minipage}{.02\textwidth} $\bu$ \linebreak \phantom{X} 
      \end{minipage}  
            & \begin{minipage}{.8\textwidth} 
              $R = \sk[x_0, \dots, x_n]\,/\,(f(x_0, \dots, x_n))$ 
              where %
              $\deg{x_i} = 1$ for all $i$ and $f$ is a homogeneous 
              polynomial of degree $d$.
              \end{minipage}  \\
       $\bu$ & $X = \Proj{R}$ is a smooth $k$-variety. \\
       $\bu$ & $\fm = (x_0, \dots, x_n)$ is the only non-regular prime of
               $R$. \\
    \end{tabular}
  \end{gathered}
\end{equation}

For each $\fp \in \Spec{R} \setminus \{\mathfrak m\}$, $R_{\fp}$ is regular.  Therefore, for
a finitely generated $R$-module $M$, the $R_{\fp}$-projective dimension of
$M_{\fp}$ is finite.  If $N$ is also finitely generated,
then for $i \gg 0$, the module $\Tor_i^R(M,N)$ is supported
on $\{\fm\}$ and hence has finite length. Thus $\theta$ is defined
for all pairs of finitely generated $R$-modules
$M$ and $N$.

The variety $X$ is smooth if and only if the radical of the
homogeneous ideal $(f, \frac{\partial f}{\partial x_0}, \dots,
\frac{\partial f}{\partial x_n})$ is $\fm$ \cite[Theorem 5.3]
{Ha}. In particular, the fourth assumption in
\eqref{assumptions} follows from the third.  Moreover,
assumptions \eqref{assumptions} remain valid upon passing to any field extension of $k$.
Further, for finitely generated $R$-modules $M$ and $N$ and any field
extension $\sk \subset \sk'$, we have 
\begin{equation} \label{theta field stable} 
  \theta^{R}(M,N) = \theta^{R \otimes_\sk \sk'}(M \otimes_\sk \sk', N \otimes_\sk \sk').
\end{equation}
In many of our results and constructions, we assume $\sk$ is
separably closed.  Some of our results apply only when $\sk = \bC$.

\subsection{Geometry}

Let $p \colon Y \to \Spec{R}$ be the blow-up of $\Spec{R}$ at the
point $\fm$, so that $Y = \Proj{\oplus_{i \geq 0} \fm^it^i}$. Note
that $\fm^i = \oplus_{j \geq i} R_j$.  The exceptional fiber of $p$ is
\begin{equation*}
  \Proj{\oplus_{i \geq 0} \fm^i/\fm^{i+1}} \iso \Proj{R} = X.
\end{equation*}
Moreover, $Y$ is the geometric line bundle over $X$ associated to the
rank one locally free coherent sheaf $\cO_X(1)$, and the inclusion $i
\colon X \into Y$ is the zero section of this line bundle
\cite[Lemma 2.2]{CHWW}. The projection $\pi \colon Y \onto X$ comes from the
inclusion of graded rings $R \into \oplus_{i\geq0} \fm^it^i$ given by
identifying $R$ with $\oplus_{i \geq 0} R_it^i$.

Assume $\sk$ is infinite, so that there is a $\sk$-rational point $Q
\in \bP^n \setminus X$.  Then linear projection away from $Q$
determines a regular map $\bP^n \setminus \{Q\} \to
\bP^{n-1}$, and we write $\rho: X \to \bP^{n-1}$ for its
restriction to $X$. The map $\rho$ is finite and dominant of degree
$d$.  The following diagrams summarize the situation:
\begin{equation} \label{spaces}  %
  \begin{gathered}
  \xymatrix{
    X 
    \ar@{^(->}@<.5ex>[r]^-{\ds{i}} 
    \ar[d]_{\ds{q}} 
    & Y 
      \ar[d]^-{\ds{p}} 
      \ar@<.5ex>@{->>}[l]^-{\ds{\pi}}   
    \\
    \Spec{\sk} \ar[r]_-{\ds{j}} 
    & \Spec{R}     
   } 
   \end{gathered}
   \qquad {\text{and}} \qquad 
   \begin{gathered}
   \xymatrix{
     X 
     \ar@{^(->}[r]
     \ar[rd]|-*+{\ds q}
     \ar[d]_-{\ds{\rho}}
     & \bP^{n} \ar[d]
     \\
     \bP^{n-1} \ar[r]_{\ds{s}} 
     & \Spec{\sk}.  
   }
   \end{gathered}
\end{equation}
The map $\rho$ will be used in the proofs of the main results of this
paper.

For a Noetherian scheme $Z$, let $\G{Z}$ denote the Grothendieck group
of coherent sheaves on $Z$.  Thus $\G{Z}$ is the abelian group
generated by isomorphism classes of coherent sheaves modulo relations
coming from short exact sequences.  We write $\K{Z}$ for the
Grothendieck group of locally free coherent sheaves on $Z$.  Recall
that $\K{Z}$ is a ring under tensor product.  If $Z$ is a
smooth $\sk$-variety, the canonical map $\K{Z} \to \G{Z}$ is an
isomorphism.  We write $\G{R}$ for $\G{\Spec{R}}$, so that $\G{R}$ is
the usual Grothendieck group of finitely generated $R$-modules.

Since $\theta$ is biadditive \cite[page 98]{Ho} and $\theta$ is
defined for all pairs of finitely generated $R$-modules, it follows
that $\theta$ determines a pairing on $G(R)$ and hence on $\G{R}_\bQ
:= \G{R} \otimes_\bZ \bQ$,
\begin{equation*}
  \theta \colon \G{R}_\bQ \tensor_\bQ \G{R}_\bQ \to \bQ.
\end{equation*}

Since $\pi: Y \to X$ is a line bundle and $X$ and $Y$ are smooth,
pull-back along $\pi$ determines an isomorphism
$
   \pi^* \colon \K{X} \map{\iso} \K{Y}, 
$
with the inverse map given by $i^*$.  The composition 
$
  \smash{\K{X} \map{i_*} \K{Y} \map{i^*} \K{X}}
$
is multiplication by the element $\alpha := [\cO_X] - [\cO_X(1)]$ of $K(X)$.
We may also describe this class as $\alpha = -[\cO_H(1)]$,
where $H \subset X$ is a general hyperplane section of $X$.
Finally, the map $ j_* \colon \G{\sk} \to \G{R} $ is torsion. Indeed,
we may find a homogeneous prime $\fp$ of height $n-1$ and a
homogeneous element $x \in \fm \setminus \fp$.  Then the short exact
sequence $\smash{R/\fp \stackrel x \into R/\fp \onto R/(x, \fp)}$
shows that $0 = [R/(x,\fp)] \in \G{R}$.  As $(x, \fp)$ is homogeneous
and $R/(x,\fp)$ has finite length, a prime filtration shows
$[R/(x,\fp)] = \length(R/(x,\fp)) \cdot [R/\fm] \in \G{R}$.  Hence the
generator $[\sk]$ for $\G{\sk}$, which maps to $[R/\fm]$, is
annihilated by $\length(R/(x,\fp)) \in \bN$.

Applying the localization sequence for $G$-theory to the left-hand
square in \eqref{spaces} yields the diagram with exact rows
$$
\xymatrix{
\ar@{-->}[r] & G_1(U) \ar[d]^{=} \ar[r] &  G_0(X) \ar[r]^{i_*} \ar[d]
& G_0(Y) \ar[d] \ar[r] & G_0(U)
\ar[r] \ar[d]^= & 0 \\
\ar@{-->}[r] & G_1(U) \ar[r] &  G_0(\Spec k) \ar[r] & G_0(\Spec R) \ar[r] & G_0(U)
\ar[r] & 0, \\
}
$$
where $U = Y \setminus X = \Spec R \setminus \Spec k$. (In this
diagram, $G_0$ is the group written as $G$ everywhere else in this paper
and $G_1$ denotes the first higher $K$-group of the abelian category
of coherent sheaves on a scheme.)
This leads to the right exact sequence
$$
G(X) \to G(Y) \oplus G(k) \to G(R) \to 0
$$
and  since $G(k)_\Q \to G(R)_\Q$ is the zero map, we obtain the right
exact sequence
$$
G(X)_\Q \map{i_*} G(Y)_\Q \map{p_*} G(R)_\Q \to 0
$$
Since $K(X) \cong G(X)$, $i^*: K(Y) \to K(X)$ is an
isomorphism (whose inverse is $\pi^*$) and $i^* \circ i_*$ is multiplication by $\alpha \in
K(X)$, we obtain the isomorphism
\begin{equation} \label{p_* pi^*}
p_* \pi^*: {\K{X}_\bQ}/{\ip{\alpha}} \map{\iso} \G{R}_\bQ,
\end{equation} 
which allows
us to regard $\theta$ as a pairing on
$\K{X}_\bQ/\ip{\alpha}$.  
One may verify that the isomorphism \eqref{p_* pi^*} is given by ``forgetting
the grading''; i.e., for a finitely generated graded $R$-module $M$
with associated coherent sheaf $\widetilde{M}$ on $X$,  
it sends $\smash{[\widetilde{M}] \in \K{X}_\bQ}$ to
$[M] \in G(R)_\bQ$.
In particular, the vector space
$G(R)_\Q$ is spanned by classes of {\em graded} $R$-modules.

\subsection{Cohomology}

Our main technique will involve factoring the $\theta$ pairing 
through cohomology (either \'etale or singular, depending on $\sk$) via
the Chern character.  In this section we review the concepts concerning these topics
that we will need. For the assertions concerning singular cohomology
made below,
we refer the reader to 
\cite{Fu} and \cite{GH}. For those concerning \'etale cohomology,
the ultimate reference is SGA \cite{SGA4V1,SGA4V2,SGA4V3,SGA5,SGA4.5},
but a good survey of this
material can be found in \cite{FK}. The features of \'etale and
singular cohomology we need are those common to any ``Weil cohomology
theory''; see \cite[\S 3]{Kl} for a precise description of what this
means.

If $\sk$ is a separably closed field, then for any prime $\ell \neq
\chr \sk$, we can consider the \'etale cohomology of $X$ with
coefficients in the $\ell$-adic rationals.  Using $\mu_r$ to denote
the \'etale sheaf of $r$-th roots of unity, the maps
$\mu_{\ell^{m+1}}^{\tensor i} \onto \mu_{\ell^m}^{\tensor i}$ given by
taking $\ell$-th powers form an inverse system of \'etale sheaves.  When $i=0$, take
$\mu_{\ell^m}^{\tensor 0}$ to be $\Z/\ell^m$ and the map
$\Z/\ell^{m+1} \onto \Z/\ell^m$ to be the canonical one.  
By definition, 
\begin{equation*}
  \H[\et]{2i}{X}{\bQl(i)} = %
  \H[\et]{2i}{X}{\bZl(i)} \tensor_{\bZl} \bQl %
  \quad \text{and} \quad %
  \H[\et]{2i}{X}{\bZl(i)} = %
  \varprojlim_m \H[\et]{2i}{X}{\mu_{\ell^m}^{\tensor i}}.
\end{equation*}
We write $\H[\et]{ev}{X}{\bQl}$ for  $\oplus_{i \geq 0} \H[\et]{2i}{X}{\bQl(i)}$, 
which is a commutative algebra under cup product $\cup$. 

When $\sk = \bC$, let $X(\bC)$ be the complex manifold associated to
$X$.  The singular cohomology $\H{\bullet}{X(\C)}{\bQ}$ of the
manifold $X(\C)$ is a graded-commutative
$\bQ$-algebra under cup
product.  We write $\H{ev}{X(\bC)}{\bQ}$ for $\oplus_{i \geq 0} \H{2i}{X(\bC)}{\bQ}$,
the even degree subalgebra of $\H{\bullet}{X(\bC)}{\bQ}$.

The Chow group of cycles modulo rational equivalence on $X$ is
$\CH[\bullet]{X}$.  Since $X$ is smooth, it is a ring under
intersection of cycles.  The \'etale and topological cycle class maps 
are ring homomorphisms 
\begin{equation*}
  \cy[\et] \colon \CH[\bullet]{X}_{\bQl} \to \H[\et]{ev}{X}{\bQl}
  \quad {\text{and}} \quad %
  \cy[\topo] \colon \CH[\bullet]{X}_\bQ \to \H[]{ev}{X(\bC)}{\bQ}
\end{equation*}
which commute with both push-forward and pull-back maps for morphisms
of smooth, projective varieties. (The map $\cy[\et]$ is defined when $k$ is separably
closed, and $\cy[\topo]$ is defined when $k = \bC$.)

Let $\ch \colon \K{X}_\bQ \to \CH[\bullet]{X}_\bQ$ be the Chern character
taking values in the Chow ring \cite[page 282]{Fu}.  The \'etale and
topological Chern characters
\begin{equation*}
  \ch[\et] \colon \K{X}_\bQ \to \H[\et]{ev}{X}{\bQl}
  \quad \text{and} \quad
  \ch[\topo] \colon \K{X}_\bQ \to \H[]{ev}{X(\C)}{\bQ}
\end{equation*}
are defined so that
\begin{equation} 
  \label{cycle class map} 
  \begin{gathered}
  \xymatrix{
   \K{X}_\bQ \ar[d]_{\ds \ch} \ar[rd]^{\ds \ch[\et]} \\
   \CH[\bullet]{X}_\bQ \ar[r]_{\ds \cy[\et]} 
   & \H[\et]{ev}{X}{\bQl}
  }
  \end{gathered}
  \quad {\text{and}} \quad 
  \begin{gathered}
  \xymatrix{
   \K{X}_\bQ \ar[d]_{\ds \ch} \ar[rd]^{\ds \ch[\topo]} \\
   \CH[\bullet]{X}_\bQ \ar[r]_{\ds \cy[\topo]} & 
   \H{ev}{X(\C)}{\bQ}
  }
  \end{gathered}
\end{equation}
commute.  These characters are ring homomorphisms from $\K{X}_\bQ$ taking
values in graded rings.  Let $\beta \in \CH[1]{X}_\bQ$ denote the class of
a generic hyperplane section of $X$.  Let $\gamma = \cy[\et](\beta)
\in H^2_\et(X, \bQl(1))$; then $\gamma$ is the \'etale cohomology
class of the divisor given by a generic hyperplane section of $X$.
Since $ch(\cO_X) = 1$ and $ch(\cO_X(1)) = e^{\beta} = 1 + \beta +
\frac{\beta}{2!} + \cdots$, 
we have
$$
ch(\alpha) = \beta \cdot u
\hskip .5in \text{  where  } \hskip .5in
u = -1 - \frac{\beta}{2!} - \frac{\beta^2}{3!} - \cdots \in CH^\bullet(X).
$$
Since $u$ is a unit in the Chow ring of $X$, the ideals
of $CH^\bullet(X)$ generated by $\ch(\alpha)$ and $\beta$ coincide.
Likewise $\ch[\topo](\alpha)$ and $\gamma$ agree up to a unit factor
in the cohomology ring of $X$.

Since $X$ is a smooth hypersurface in $\bP^n$, Poincar\'e
duality and the weak Lefschetz theorem show that the even degree
\'etale cohomology groups of $X$ are spanned by powers of $\gamma$,
except possibly in degree ${n-1}$ when $n$ is odd.  That is,
the following equations hold: 
\begin{equation}\label{etale cohomology}
  \H[\et]{2i}{X}{\bQl(i)} = \bQl \cdot \gamma^i, %
  \quad %
  \text{for all $i$ except when $n$ is odd and $2i=n-1$.}
\end{equation}
When $\sk = \bC$, we also write $\gamma$ for the element
$cy_{\topo}(\beta) \in \H{2}{X(\C)}{\bQ}$.  Equations analogous to
\eqnref{etale cohomology} hold for singular cohomology when $k = \bC$.
Since $\ch(\alpha)$ coincides with $\beta$ up to a unit factor, 
there are
induced maps on quotient rings from ${\K{X}_\bQ/\ip{\alpha}}$ to each
of $\CH[\bullet]{X}_\bQ/\ip{\beta}$, $\H[\et]{ev}{X}{\bQl}/\ip{\gamma}$,
and $\H{ev}{X(\C)}{\bQ} /\ip{\gamma}$.  We write these induced maps
also as $\ch, \ch[\et]$, and $\ch[top]$, respectively.

Recall that for a (possibly singular) variety $Y$, the Grothendieck-Riemann-Roch isomorphism 
$$
\tau: G(Y)_{\bQ} \map{\iso} CH^\bu(Y)_\Q,
$$
is functorial for pushforwards along proper maps
\cite[Corollary 18.3.2]{Fu}.
If $Y$ happens to be smooth, so that $K(Y) \cong G(Y)$, we write $\tau$
also for the composition of isomorphisms
$$
K(Y)_{\bQ} \map{\iso} G(Y)_{\bQ} \map{\tau} CH^\bu(Y)_\Q.
$$
It is useful for our purposes to compare $\ch$ with $\tau$ for the
variety $X$.

\begin{lem} \label{ch and tau agree} %
  For $X$ as in \eqref{assumptions}, the isomorphisms $\ch$ and $\tau$ from
  $\K{X}_\bQ$ to $CH(X)_\bQ$ each map $\ip{\alpha}$ isomorphically
  onto $\ip{\beta}$. Moreover, they induce the same isomorphism
$
  {\K{X}_\bQ}/{\ip{\alpha}} \stackrel{\iso}{\lra}
  {\CH[\bullet]{X}_\bQ}/{\ip{\beta}}.
$
\end{lem}

\begin{proof}
Since $X$ is smooth, $\tau(a) = \ch(a) \, \cup\, \Td{\cT_X}$ \cite[page 287]{Fu} where
$\Td{\cT_X} \in \CH[\bullet]{X}$ is the Todd class of the tangent
bundle $\cT_X$ of $X$ \cite[Example 3.2.13]{Fu}.  Since $X$ is a smooth hypersurface of degree
$d$ in $\bP^{n-1}$, we have $[\cT_X] = (n+1)[\cO_X(1)] - [\cO_X] -
[\cO_X(d)] \in \K{X}$.  (See \cite[Examples 3.2.11, 3.2.12 \& Appendix B.7.1]{Fu}.)  Hence $\Td{\cT_X} =
\Td{\cO_X(1)}^{n+1}/\Td{\cO_X(d)}$.  Since $\Td{\cO_X(1)}$ and
$\Td{\cO_X(d)}$ are in $1 + \beta \CH[\bullet]{X}_\bQ$, so is $ \Td{\cT_X}$.
\end{proof}

For a smooth projective variety $Z$ over a separably closed field
$\sk$ with structure map $q \colon Z \to \Spec{\sk}$, we write 
\begin{equation*}
  \sint_Z \colon \H[\et]{ev}{Z}{\bQl} \to \bQl
\end{equation*}
for push-forward along $q$; it takes values in
$\H[\et]{ev}{\Spec{\sk}}{\bQl} = \H[\et]{0}{\Spec{\sk}}{\bQl} = \bQl$.
For the analogous map on singular cohomology, we also write $\sint_Z$.

\begin{example}\label{integration and gamma} 
  As $\gamma$ is the class of a hyperplane on the $n-1$ dimensional
  variety $X$, the degree of $X$ is equal to $\int_X \gamma^{n-1} =
  d$.  
\end{example}

\section{Statement of main results} \label{statement section}

We continue with the notation from Section \ref{background}. In this
section, 
we state the main results of the paper, postponing several of the proofs until Section \ref{proof section}.

In light of the isomorphism
\eqref{p_* pi^*}, we write simply $\theta(x,y)$ for
$\theta(p_*\pi^*x,p_*\pi^*y)$ when $x, y \in \K{X}_\bQ/\ip{\alpha}$.
The following proposition shows that the pairing $\theta(x,y)$ on
$G(R)_\bQ$ factors through cohomology.
\begin{prop}\label{thm} %
  Let $R$ and $X$ be as in \eqref{assumptions} with $\sk$ a separably
  closed field.  For any $x$ and $y$ in $\K{X}_\bQ/\ip{\alpha}$, we
  have
 \begin{equation*}
  \theta(x, y) = \int_{\bP^{n-1}} \Big( \rho_*(\ch[\et]{x})
    \cup\rho_*(\ch[\et]{y}) - d \cdot \rho_*\big(\ch[\et]{\big(x\cdot y)}\big) \Big).
 \end{equation*}

If $k = \C$, the analogous formula involving $ch_{\topo}$ also holds. 
\end{prop}

This proposition is useful in proving the following theorem, which establishes Dao's \conjref{Dao Conj}
for those isolated hypersurface singularities satisfying \eqref{assumptions}.  

\begin{thm}\label{n even result} %
  Let $R$ and $X$ be as in \eqref{assumptions} with $\sk$ an arbitrary field.
  If $n$ is even, then $\theta$ vanishes; i.e., for every pair of
  finitely generated modules $M$ and $N$, 
\begin{equation*}
  \theta(M,N) = \length(\Tor_{2i}^R(M,N)) -
  \length(\Tor_{2i+1}^R(M,N)) = 0 \qquad {\text{for all $i \gg 0$}}.  
\end{equation*}
\end{thm}

We now give a precise description of $\theta$ when $n$ is odd and
$\sk$ is a separably closed field.  In this case, we define a symmetric
pairing $\theta_\et$ on the $\bQl$-vector space
$\H[\et]{n-1}{X}{\bQl(\frac{n-1}{2})}$ by setting 
\begin{equation} \label{theta on cohomology} %
 \theta_\et(a, b) = \left(\sint_X a \cup \gamma^{\frac{n-1}{2}}\right) %
  \left(\sint_X b \cup \gamma^{\frac{n-1}{2}}\right) %
  - d \left(\sint_X a \cup b \rule{0em}{1em} \right).
\end{equation}
When $n = 1$, by $\gamma^0$ we mean $1 \in \H[\et]{0}{X}{\bQl}$.  When
$\sk = \bC$, the same expression on the right-hand side of (\ref{theta on cohomology}) defines a
symmetric bilinear pairing $\theta_\topo$ on the singular cohomology
group $\H{n-1}{X(\bC)}{\bQ}$.  By Example \ref{integration and gamma},
if either $a$ or $b$ is a $\bQl$-multiple of $\gamma^{\frac{n-1}{2}}$,
then $\theta_\et(a,b) = 0$; similarly, when $\sk = \bC$, if $a$ or $b$
is a $\bQ$-multiple of $\gamma^{\frac{n-1}{2}}$, then
$\theta_\topo(a,b) = 0$.  Thus $\theta_\et$ and $\theta_\topo$
induce pairings on
\begin{equation*}
  \frac{\H[\et]{n-1}{X}{\bQl(\frac{n-1}{2})}}{\bQl \cdot
    \gamma^{\frac{n-1}{2}}} %
  \qquad {\text{and}} \qquad %
  \frac{\H{n-1}{X(\bC)}{\bQ}}{\bQ \cdot
    \gamma^{\frac{n-1}{2}}}, \quad {\text{ respectively;}}
\end{equation*}
we retain the $\theta_\et$ and $\theta_\topo$ notation for these
pairings.  

\begin{thm}\label{n odd thm} %
Let $R$ and $X$ be as in \eqref{assumptions} with $\sk$ a separably closed
field. If $n$ is odd, then there is a commutative diagram: 
\begin{equation*}
\begin{gathered}
\xymatrix{
\G{R}_\bQ^{\tensor 2} \ar[d]_-{\ds\theta}
  & \displaystyle{\left(\frac{\K{X}_\bQ}
                             {\ip{\alpha}}\right)^{\tensor 2}} 
    \ar[l]_-{\ds (p_*\pi^*)^{\tensor 2}}^-{\ds\iso}
    \ar[d]^-{\ds \left(ch_{\et}^{\frac{n-1}{2}} \right)^{\tensor 2}} 
 \\
\bQl
     & \displaystyle{\left(\frac
                         {\H[\et]{n-1}{X}{\bQl(\frac{n-1}{2})}} 
                         {\bQl \cdot \gamma^{\frac{n-1}{2}}}\right)^{\tensor 2}}
       \ar[l]^-{\ds \theta_\et} 
}
\end{gathered}
\end{equation*}

When $\sk = \C$, the analogous diagram involving $H^{n-1}(X(\C), \Q)$,
$ch_{\topo}^{\frac{n-1}{2}}$, and $\theta_\topo$  also commutes; i.e., $\theta \circ (p_* \pi^*)^{\tensor 2} = %
\theta_\topo \circ (\ch[\topo]^{\frac{n-1}{2}})^{\tensor 2}$.    
\end{thm}

When $\sk = \C$, we use the Hodge-Riemann bilinear relations \cite[page 165]{Vo} to
analyze the pairing $\theta_\topo$ in more detail, obtaining:  

\begin{thm}\label{n odd pos defn} %
  For $R$ and $X$ as in \eqref{assumptions} with $\sk = \C$, if $n$ is
  odd, then the restriction of the pairing
  $(-1)^{\frac{n+1}{2}}\theta_\topo$ to
\begin{equation*}
\im \left(\ch[\topo]^{\frac{n-1}{2}} \colon 
\frac{\K{X}_\bQ}{\ip{\alpha}} 
\lra 
\frac{\H{n-1}{X(\bC)}{\bQ}}{\bQ \cdot \gamma^{\frac{n-1}{2}}}
\right)
\end{equation*}
is positive definite; 
i.e., for $v$ in this image,
$(-1)^{\frac{n+1}{2}}\theta_\topo(v,v) \geq 0$ with equality holding
if and only if $v = 0$.

In particular, for $R$ as in \eqref{assumptions} with $\sk$ an arbitrary field
of characteristic zero, the pairing $(-1)^{\frac{n+1}{2}}\theta$ on
$\G{R}_\bQ$ is positive semi-definite.
\end{thm}

\begin{rem} If the ``Hodge standard conjecture'' 
for $\ell$-adic \'etale cohomology were known (see \cite[\S 5]{Kl}),
then the evident analogue of this theorem involving an algebraically closed
field of characteristic $p$ could
be proven. The proof would be nearly identical to the one given in
Section \ref{proof section} below.
\end{rem}
 
With Theorem \ref{n odd pos defn} serving as evidence, we propose the following conjecture
about the $\theta$ pairing in general. 

\begin{conj} \label{Conj819} %
  Let $S$ be an admissible isolated hypersurface singularity of
  dimension $n$. If $n$ is odd, then $(-1)^{\frac{n+1}{2}}\theta$ is
  positive semi-definite on $\G{S}_\bQ$.
\end{conj}

The results above can be applied to give a relation between the theta pairing and Chow groups.

\begin{cor} \label{cor theta on Chow X} %
  Let $R$ and $X$ be as in \eqref{assumptions} with $k$ an arbitrary field,
  and assume $n$ is odd. Then the theta pairing is induced from a
  pairing, which we call $\theta_{CH(X)}$, on $\left(\frac{\CH[\bullet]{X}_\bQ}
                  {\ip{\beta}} \right)^{\frac{n-1}{2}}$, the $\frac{n-1}{2}$ part of the graded ring $\CH[\bullet]{X}_\bQ/\ip{\beta}$.  That is,
  there is a commutative diagram:
  \begin{equation} \label{theta on Chow X} %
 \begin{gathered}
  \xymatrix{
  G(R)_\bQ^{\tensor 2} 
  \ar[d]_-{\ds \theta} 
  & \left( {\ds \frac{\K{X}_\bQ}{\ip{\alpha}}} \right)^{\otimes 2} 
    \ar[l]_-{\ds (p_* \pi^*)^{\tensor 2}}^-\iso
    \ar@{->>}[d]^-{\ds \left(\ch^{\frac{n-1}{2}}\right)^{\otimes 2}}
  \\
  \bQ
  & \left( \left({\ds \frac{\CH[\bullet]{X}_\bQ}
                  {\ip{\beta}} 
                  }\right)^{\frac{n-1}{2}} \right)^{\tensor 2} 
     \ar[l]^-{\ds \theta_{CH(X)}} 
  }
 \end{gathered}
\end{equation} 
\end{cor}

\begin{proof}
We first argue under the assumption that $\sk$ is separably closed.
In this case, define
$ \theta_{CH(X)} $ to be $\theta_\et \circ
\left(cy_\et^{\frac{n-1}{2}}\right)^{\otimes 2}$.
By Theorem \ref{n odd thm} and \eqref{cycle class map},
the diagram similar to \eqref{theta on Chow X}, but with the
$\bQ$ in the lower-left corner replaced by $\bQl$, commutes.  As
the image of $\theta \circ (p_*  \pi^*)^{\tensor 2}$ is
contained in $\bQ \subseteq \bQl$, so too is the image of $\theta_{CH(X)}$.

For an arbitrary field $\sk$, let $\sk_{sep}$ be a separable closure of
$\sk$, let $X_{sep} = X \otimes_{\sk} \sk_{sep}$, and let $\theta_{CH(X)}$ be 

\begin{equation*} 
  \xymatrix{ \ds
\left( \left({\ds \frac{\CH[\bullet]{X}_\bQ}
                  {\ip{\beta}} 
                  }\right)^{\frac{n-1}{2}} \right)^{\tensor 2} 
     \ar[r] & \left( \left({\ds \frac{\CH[\bullet]{X_{sep}}_\bQ}
                  {\ip{\beta}} 
                  }\right)^{\frac{n-1}{2}} \right)^{\tensor 2} 
                   \ar[rr]^-{\ds \theta_{CH(X_{sep})}}
    & & \bQ.
  }
\end{equation*}

The commutativity of \eqref{theta on Chow X} follows from the fact
that $p_* \pi^*$ and $ch$ are natural with respect to pull-back
along $X_{sep} \to X$.
\end{proof}

Since $p_* \pi^*$ induces an isomorphism from 
$\left({\CH[\bullet]{X}_\bQ}/{\langle \beta \rangle}\right)^{\frac{n-1}{2}}$ to
$\CH[\frac{n-1}{2}]{R}_\bQ$, the previous Corollary shows that the
theta pairing on $G(R)_\bQ$ factors through a pairing on
$\CH[\frac{n-1}{2}]{R}_\bQ$.  

\begin{cor} \label{theta on Chow R} %
  Let $R$ be as in \eqref{assumptions} with $k$ an arbitrary field. If $n$ is
  odd, then there exists a pairing $\theta_{CH(R)}$ on
  $CH^{\frac{n-1}{2}}(R)_\bQ$, which corresponds to the pairing $\theta_{CH(X)}$ under the isomorphism $p_*  \pi^*$, such that the triangle
   \begin{equation*}
 \xymatrix{ 
 \G{R}_\bQ^{\tensor 2} 
 \ar@{->>}[rr]^-{\ds(\tau^{\frac{n-1}{2}})^{\tensor 2} } 
 \ar[rd]_(.45){\ds\theta} 
   & & \CH[\frac{n-1}{2}]{R}_\bQ^{\tensor 2} 
       \ar[ld]^(.4){\ds\theta_{CH(R)}} \\
   & \bQ
 }
 \end{equation*}
 commutes. 
Here, $\tau^{\frac{n-1}{2}}$ is the degree $\smash{\frac{n-1}{2}}$
 component of the Grothendieck-Riemann-Roch isomorphism $\smash{\tau \colon
 \G{R}_\bQ \map{\iso} \CH[\bullet]{R}_\bQ}$. 

In particular, if $\CH[\frac{n-1}{2}]{R}$
 is torsion, then $\theta = 0$.
 \end{cor}

 \begin{proof}    
   We claim that the diagram
\begin{equation} \label{E999} %
 \begin{gathered}
 \xymatrix{
 {\K{X}_\bQ}/\ip{\alpha} 
 \ar[rr]_-\iso^-{\ds p_* \pi^*} 
 \ar[d]^{\ds ch}_\iso 
 & & \G{R}_\bQ 
     \ar[d]^{\ds \tau}_\iso \\
 \CH[\bullet]{X}_\bQ/\ip{\beta}
 \ar[rr]^\iso_-{\ds p_* \pi^*} && \CH[\bullet]{R}_\bQ
 }
 \end{gathered}
\end{equation}
commutes. Granting this, the result follows from
\corref{cor theta on Chow X}, since the bottom arrow in this diagram
is graded.  In more detail,  it follows from the equalities
$\theta = \theta_{CH(X)} \circ ch^{\frac{n-1}{2}} \circ (p_*\pi^*)^{-1} =  \theta_{CH(X)} \circ
 (p_*\pi^*)^{-1} \circ \tau^{\frac{n-1}{2}} = \theta_{CH(R)} \circ \tau^{\frac{n-1}{2}}$ (suppressing the $^{\tensor 2}$ notation).
 
 To show that (\ref{E999}) commutes, we consider the diagram below, where the
 left-hand square commutes by
   \cite[Theorem 18.2(3)]{Fu} and the right-hand one commutes by \cite[Theorem 18.2(1)]{Fu}.
 \begin{equation} \label{Fulton squares}
 \begin{gathered}
 \xymatrix{
 \K{X}_\bQ \ar[rr]^-{\ds\pi^*}_\iso \ar[d]^{\ds\tau}_{\iso} 
 & & \K{Y}_\bQ \ar[r]^-{\ds p_*} \ar[d]^-{\ds\tau}_\iso
     & \G{R}_\bQ \ar[d]^-{\ds\tau}_\iso \ar[d] \\
 \CH[\bullet]{X}_\bQ \ar[rr]^-{\ds\Td{\cT_{\pi}}\pi^*}_\iso 
 & & \CH[\bullet]{Y}_\bQ \ar[r]^-{\ds p_*} 
     & \CH[\bullet]{R}_\bQ
 }
 \end{gathered}
 \end{equation}
 Here, $\cT_\pi$ is the relative tangent bundle of $\pi$ and
 $\Td{\cT_\pi}$ is its Todd class \cite[Example 3.2.4]{Fu}, which is a unit in the ring
 $\CH[\bullet]{Y}_\bQ$. Since $\pi$ is the geometric line bundle
 associated to the invertible sheaf $\cO_X(1)$, it follows that
 $\cT_\pi \iso \pi^* \cO_X(-1)$ and hence that $\Td{\cT_\pi} = \pi^*
 \Td{\cO_X(-1)}$. In particular, the map $\Td{\cT_\pi} \pi^*$ in
\eqref{Fulton squares} is 
\begin{equation*}
 \pi^*\left(\frac{-\beta}{1 - \exp(\beta)}\right)\pi^* = 
  \left(1 - \frac12 (\pi^*\beta) + \frac{1}{12} (\pi^*\beta)^2 
  - \cdots\right)\pi^*,
\end{equation*}
and so upon modding out by the ideals generated by
$\beta$ and $\Td{\cT_\pi}\pi^*(\beta)$, the map $\Td{\cT_\pi} \pi^*$ coincides
with $\pi^*$.  By Lemma \ref{ch and tau agree}, the
left-most vertical map $\tau$ in \eqref{Fulton squares} sends
$\ip{\alpha}$ isomorphically onto $\ip{\beta}$ and the maps $\tau$ and
$ch$ coincide as maps on the quotients. The commutativity of
\eqref{E999} follows.
\end{proof}

\begin{rem} To provide clarity for the relations among the above results,
we summarize them in the diagram below.  
For the sake of simplicity, the notation $^{\tensor 2}$ is suppressed.
Note that all arrows to $\bQ$ represent interpretations of $\theta$.  Recall that the arrows $\ch[\et]$ and
$\ch[\topo]$ require assumptions on $\sk$, namely that it is separably closed and equal to
$\bC$, respectively.
\end{rem}

\begin{equation}\label{all thetas}
\begin{gathered}
\xymatrix{
\dfrac{\K{X}_\bQ}{\ip{\alpha}}
\ar[drrr(.95)]^-{\theta} 
\ar[rr]^-{\iso} 
\ar[dd]_{\chj[]{\frac{n-1}{2}}} 
    & &\G{R}_\bQ 
        \ar[r]^-{\tau^{\frac{n-1}{2}}} 
        \ar[dr(.93)]^-\theta        
        & \CH[\frac{n-1}{2}]{R}_\bQ 
            \ar[d]^{\theta_{\CH[]{R}}}        
        \\ 
   & & \dfrac{\H[\et]{\frac{n-1}{2}}{X}{\bQl}}{\bQl\gamma^{\frac{n-1}{2}}} 
       \ar[r]|-{\theta_{\et}}      
      & \bQ \\
\left(\dfrac{\CH[\bullet]{X}_\bQ}{\ip{\beta}}\right)^{\frac{n-1}{2}} 
\ar[rrru(.94)]_-{\theta_{CH(X)}} 
\ar[rr]_-{\cy[\topo]^{\frac{n-1}{2}}} 
\ar[urr(.8)]^-{\cy[\et]^{\frac{n-1}{2}}}
  & & \dfrac{\H{\frac{n-1}{2}}{X(\C)}{\Q}}{\Q\gamma^{\frac{n-1}{2}}} 
      \ar[ur(.87)]_-{\theta_{\topo}} 
     \\
}
\end{gathered}
\end{equation}

We end this section with some further applications and examples.  Recall that
homological equivalence of cycles is defined 
so that the image of the cycle class map
$\cy[\topo]$ is isomorphic to the group of cycles modulo homological
equivalence: $\im \cy[\topo] \iso \CH[\bullet]{X}/({\text{hom} \sim
  0})$; see \cite[Definition 19.1]{Fu}.
 
\begin{cor} \label{n odd theta = 0} Let $R$ and $X$ be as in
  \eqref{assumptions} with $k = \C$. The rational vector space
$\CH[\frac{n-1}{2}]{X}_\bQ/({\text{hom} \sim 0})$ 
is spanned by the $\frac{n-1}{2}$-st
multiple of the class of a hyperplane section if and only if $\theta =
0$.  

For $n=3$, the divisor class group $\CH[1]{R}$ is
torsion if and only if $\theta = 0$.
\end{cor}

\begin{proof}
The first assertion follows from Theorem \ref{n odd pos defn} and
Corollary \ref{cor theta on Chow X}.  Since $X$ is a smooth hypersurface in projective space, homological
  and rational equivalence coincide on codimension one cycles: $
  \CH[1]{X}_\bQ $ is isomorphic to $ \CH[1]{X}_\bQ/({\text{hom} \sim
    0})$.  The second assertion therefore follows from the isomorphism
  $\CH[1]{X}_\bQ/\bQ \cdot \beta \iso \CH[1]{R}_\bQ$. \qedhere 
\end{proof}

\begin{example}
  If $n=1$, then since $X$ is smooth, it consists of $d$ distinct
  points $P_1, \ldots, P_d$ in $\bP^1$.  We have $\K{X}_\bQ =
  \CH[\bullet]{X}_\bQ = \CH[0]{X}_\bQ = \bQ^d$ and the map $\ch[\et]
  \colon \CH[\bullet]{X}_{\bQl} \to \H[\et]{ev}{X}{\bQl} =
  \H[\et]{0}{X}{\bQl}$ is an isomorphism.  A basis of $CH^0(X)_\bQ$
  is given by the classes of the $P_i$'s.  Theorem \ref{n odd thm}, or
  direct calculation, gives
  \begin{equation*}
    \theta(P_i, P_j) = 1 - d \cdot \delta_{ij}.
  \end{equation*}
  Since $\beta^0 = 1 = P_1 + \cdots + P_d$, a basis for
  $\CH[\bullet]{X}_\bQ/\ip{\beta^0}$ is given by $(P_1 -
    P_i)/\sqrt{2d-2}$ for $i > 1$.  With this basis, $\theta$ is
  represented by the matrix $- I_{d-1}$, and so is negative definite.
\end{example}

\begin{example} \label{Ex817} %
  Let $\sk$ be a separably closed field and $n=3$, so that $X$ is a
  smooth surface of degree $d$ in $\bP^3$.  In this case, it follows from
  Corollary \ref{cor theta on Chow X} that the $\theta$ pairing on
  $\K{X}_\bQ$ is induced from a pairing $\theta_{CH(X)}$ on
  $\Pic(X)_\bQ = \CH[1]{X}_\bQ$ via the map
\begin{equation*}
  \xymatrix{
    \K{X}_\bQ 
    \ar[rr]^-{\ds \ch^1 = c_1} 
    & & \CH[1]{X}_\bQ,
  } 
\end{equation*}
where $c_1$ is the first Chern class.
Observe that for a curve $C$ on $X$, we have $\sint_X cy_\et(C) \cup
\gamma^{\frac{n-1}{2}} = \deg{C}$.  Thus, from \thmref{n odd thm}
and \corref{cor theta on Chow X}, we get
\begin{equation} \label{E817b}
d \cdot (C \cap D) = \deg(C) \deg(D) - \theta_{CH(X)}(C,D),
\end{equation}
where  $C$ and $D$ are curves on $X$ and $C \cap D$ denotes the number
of points of their intersection, counted with multiplicity. 

It is perhaps useful to think of \eqref{E817b} as a generalization of
B\'ezout's Theorem.  If $d =1$, then $X = \bP^2_\C$ and $\theta \equiv
0$ (since $R$ is regular), so that the resulting equation is the
classical B\'ezout's Theorem. For $d > 1$, $\theta$ gives the
``error'' term of this generalized version of B\'ezout's Theorem.

We also remark that the positive definiteness of 
\begin{equation*}
  \theta_{CH(X)}(C,D) = 
  \deg(C) \deg(D) - d \cdot (C \cap D) 
\end{equation*}
on $\CH[1]{X}_\bQ$ is a consequence of the Hodge Index
Theorem; see \cite[page 165]{Vo}.  
\end{example}

\begin{example}  In these examples, lower case letters refer to
  the images of the upper case variables in the quotients.  
\begin{enumerate} 
\item $(n = 3)$: Let $R$ be the ring $\bC[X,Y,U,V]/(XU + YV)$, so that
  $X = \bP^1 \times \bP^1$, embedded in $\bP^3$ via the Segre
  embedding. Since $\CH[1]{R} \iso \bZ$, $\theta$ does not vanish, and
  moreover it is positive definite since $\frac{n+1}{2}$ is even.  Set
  $M = R/(x, y)$.  Matrices for the minimal resolution of $M$
  eventually alternate between $\left[\begin{smallmatrix}
      x & y \\
      v & -u
   \end{smallmatrix}\right]$ 
  and
  $\left[\begin{smallmatrix} 
    u & y \\
   v & -x
   \end{smallmatrix}\right]$.  
 It is now easy to calculate that $\theta(M,M) = 1$.

\item $(n = 5)$: Let $R$ be $\bC[X, Y, Z, U, V, W]/ (XU
  + YV + ZW)$.  Since $\frac{n+1}{2}$ is odd, $\theta$ is
  negative definite. Set $M = R/(x, y, z)$.  Then $\theta(M,M)
  = -1$.  Here the matrices for the minimal resolution of $M$
  eventually alternate between 
\begin{equation*}
 \left[\begin{smallmatrix}
  x & y &z & 0 \\
  v &-u & 0 & z \\
  -w & 0 & u & y \\
  0 & -w & v & -x 
 \end{smallmatrix}\right]
\quad {\text{and}} \quad 
 \left[\begin{smallmatrix}
  u & y & -z & 0 \\
  v &- x & 0 & -z \\
 w & 0 & x & y \\
  0 & w & v & -u
 \end{smallmatrix}\right].
\end{equation*}
\end{enumerate}
\end{example}

\begin{example} \label{Buchweitz ex} %
In unpublished notes \cite[\S10.4]{Bu}, Buchweitz studies the
  Herbrand difference pairing $h(-,-)$ on $\G{R}$ for $R$ and $X$ as
  in \eqref{assumptions} with $n = 3$ and $f(x_0, x_1, x_2, x_3)$ of
  degree three. As with the pairing $\theta$, $h$ can be
  interpreted as a pairing on $\K{X}$ induced from a pairing
  on $\Pic(X) = CH^1(X)$.  In this example, $\Pic(X)$ is free of rank
  six, with basis given by the classes of six lines on $X$.  Let $M$
  and $M'$ be cyclic modules defining two lines $L$ and $L'$ on $X$.
  With $\omega$ as the canonical divisor of $X$, by direct
  calculation, Buchweitz obtains
\begin{equation*}
  h(M, M') = - \frac{1}{3} (3L + \omega) \cap (3L' + \omega),
\end{equation*} 
and remarks that ``there should exist a more conceptual proof for
this''. We show how our Theorem \ref{n odd thm} leads to Buchweitz's formula.

Since $X$ is a degree three hypersurface in $\bP^3$, we have 
$\omega = -\beta$ and hence  $L \cap \omega = L' \cap \omega = -1$
and $\omega \cap \omega= 3$.  It follows that Buchweitz's formula for $h$ 
is equivalent to
\begin{equation*}
  h(M,M') 
  = 1 - 3(L \cap L').
\end{equation*}
As mentioned in the introduction, $\theta$ and $h$ coincide on
$CH^1(X)$, and hence they coincide as pairings on $G(R)$ in this example. Buchweitz's formula for $h$ thus
follows from the formula for $\theta$ given in Example \ref{Ex817} (with $C = L$ and $D = L'$, so that $\deg(C) =
\deg(D) = 1$ and $d=3$).
\end{example}

The following corollary shows that, at least when $\chr(k) = 0$, to
check rigidity of a module one needs only to check $\theta$ of the
module against itself.

\begin{cor} \label{n odd theta rigidity} %
  Let $R$ be as in \eqref{assumptions} with $\sk$ of characteristic $0$ and
  let $n$ be odd.  If $M$ is a finitely-generated $R$-module with
  $\theta(M,M) = 0$, then $M$ is rigid.
\end{cor}

\begin{proof} 
  If $\phi$ is a positive semi-definite form on a $\bQ$-vector space
  $V$ and $v \in V$, then $\phi(v,v) = 0$ implies $\phi(v,-) \equiv
  0$.  Since $(-1)^{\frac{n+1}{2}}\theta$ is positive semi-definite by
  Theorem \ref{n odd pos defn}, if $\theta(M, M) = 0$, then
  $\theta(M,N) = 0$ for all finitely generated $R$-modules $N$.  The
  result now follows from \cite[Proposition 2.8]{Dao1}.
\end{proof}

\begin{rem}  
  \thmref{n even result} and \cite[Proposition 2.8]{Dao1} imply that when
  $R$ is as in \eqref{assumptions}, with $\sk$ arbitrary and $n$ even, then
  every finitely generated $R$-module is rigid.  
\end{rem}

\section{Hilbert series, Hilbert polynomials and related invariants} 
\label{Hilbert section}

This section establishes some technical results that will be used in Section \ref{proof section} to 
prove the main results of this paper.  Let $R$ be as in \eqref{assumptions}. For the first part of this
section, we assume only that the field $\sk$ is infinite. Recall that
the map $\rho \colon X \to \bP^{n-1}$ is defined as projection away from a
$k$-rational point of $\bP^n \setminus X$.

Throughout this section we identify $K(\bP^{n-1})_\bQ$ with
$\bQ[t]/(1-t)^n$ 
under the ring isomorphism
sending $t$ to $[\cO_{\bP^{n-1}}(-1)]$
(see \cite[Exercise III.5.4]{Ha}).  For example, given $x \in \K{X}_\bQ$, we interpret $\rho_*(x)$ as
being a truncated polynomial in $t$, i.e., an element of
$\bQ[t]/(1-t)^n$, and we use the fact that the ring homomorphism $\rho^*$
satisfies $\rho^*(t) = [\cO_X(-1)] \in \K{X}_\bQ$.

For a finitely generated graded $R$-module $M$, its {\emph{Hilbert series}}
\begin{equation*}
  \Hs[M]{t} = \sum_{l \in \Z} \dim[\sk]{M_l} t^l
\end{equation*}
is a rational function with a pole of order equal to $m := \dim{M}$ at
$t=1$. In fact, 
\begin{equation*}
  \Hs[M]{t} = \frac{e_M(t)}{(1-t)^{m}}, 
\end{equation*}
where $e_M(t)$ is a Laurent polynomial \cite[(1.1)]{AB}.  
The {\emph{Hilbert polynomial}} of $M$ is the 
polynomial $\Hp[M]{l}$ of degree $m-1$ such that
\begin{equation*}
  \Hs[M]{t} = \text{ some Laurent polynomial in $t$ } 
              + \sum_{l \geq 0} \Hp[M]{l}t^l.
\end{equation*}
For $j \geq 1$, let $q_j(l)$ be the degree $j-1$ polynomial with $\bQ$
coefficients given by
\begin{equation*}
q_j(l) = {l + j - 1 \choose l} = \frac{(l+j-1) \cdots (l+1)}{(j-1)!}.  
\end{equation*}
For $j \leq 0$, let $q_j = 0$.  So $q_j(0) = 1$ for all $j \geq 1$,
and $q_j(0) = 0$ for all $j \leq 0$.  Thus 
\begin{equation} \label{E814b}
\frac{1}{(1-t)^j} = \sum_{l \geq 0} q_j(l) t^l.  
\end{equation}

We now assume $M$ is non-negatively graded (i.e., $M_l = 0$ for all $l
< 0$).  Then $e_M(t) = \sum_{i \geq 0} a_i (1 - t)^i$ is a polynomial.
Hence for some rational numbers $a_i$,
\begin{equation*}
\Hs[M]{t} =
 \frac{a_0}{(1-t)^{m}} +
\frac{a_1}{(1-t)^{m-1}} + \cdots +
\frac{a_{m-1}}{(1-t)^{1}} + \text{ a polynomial}.
\end{equation*}
Using \eqref{E814b}, we get
\begin{equation} \label{E10}
\Hp[M]{l} = a_0 q_{m}(l) + \cdots + a_{m-1} q_1(l).
\end{equation}

Recall that the {\em first difference} of a polynomial $q(l)$ is the
polynomial $q^{(1)}(l) = q(l) - q(l-1)$, and recursively one defines $ q^{(i)} =
(q^{(i-1)})^{(1)}$.  For all $j$, by induction on $i$, one may prove 
$
  q_j^{(i)}(l) = q_{j-i}(l),
$
and so, 
\begin{equation*}
  q_j^{(i)}(0) = q_{j - i}(0) = 
\begin{cases}
  1 & \text{if $j > i$ and} \\
  0 & \text{if $j \leq i$}.
\end{cases}
\end{equation*}
Thus, still assuming $M$ is non-negatively graded, \eqref{E10} gives 
\begin{equation} \label{E11} %
\Hpi[M]{i}{0} = a_0 + \cdots + a_{m -i - 1}.
\end{equation}

The lemma below shows that the coefficients of $\rho_*([\widetilde{M}])$ in the basis $1$, $1-t$,
$\cdots$, $(1-t)^{n-1}$ of $\bQ[t]/(1 - t)^n \iso \K{\bP^{n-1}}_\bQ$
coincide, up to order, with the coefficients of $\Hp[M]{l}$ in the
basis $q_1(l), \dots, q_n(l)$ of $\bQ$-polynomials of degree at most
$n - 1$; see \eqref{E12}.

\begin{lem} \label{rho and HS} %
Let $R$ be as in \eqref{assumptions} with $\sk$ infinite,
and let $M$ be a finitely generated graded $R$-module.  Then 
\begin{equation*}
\rho_*([\widetilde{M}]) = (1-t)^n\Hs[M]{t}
\quad \text{in $\K{\bP^{n-1}}_\bQ \iso \bQ[t]/(1-t)^n$.}
\end{equation*}
In particular,
\begin{equation*}
\rho_*(1)=
\rho_*([\cO_X])) = e_R(t) = 1 + t + \cdots +
  t^{d-1}
\quad \text{in $\bQ[t]/(1-t)^n$.}
\end{equation*}
\end{lem}

\begin{proof} We leave the second remark to the reader.  

  Using $\rho^*(t) = [\cO_X(-1)]$ and the projection formula, we get
\begin{equation*}
  \rho_*([\widetilde{M}(-i)]) 
  = \rho_*([\widetilde{M}] \cdot \rho^*(t^i)) 
  = \rho_*([\widetilde{M}]) t^i
\quad \text{in $\bQ[t]/(1-t)^n$.}
\end{equation*}
Likewise,
\begin{equation*}
\Hs[M(-i)]{t}(1-t)^n = t^i \Hs[M]{t}(1-t)^n
\quad \text{in $\bQ[t]/(1-t)^n$.}
\end{equation*}
It follows that the lemma holds for $M$ provided it holds for $M(-i)$
for any $i$.  So taking $i$ sufficiently large, we may assume
$M$ is non-negatively graded.
Therefore, we have \eqnref{E11}.  Further
\begin{align*}
\Hs[M]{t}(1-t)^n = & ~e_M(t)(1-t)^{n-m} \\
                 = & ~a_0(1-t)^{n-m} + a_1(1-t)^{n-m+1} 
                     + \cdots + a_{m-1}(1-t)^{n-1} \\
                   & + \text{higher order terms in\ } (1 - t).
\end{align*}
There are rational numbers $b_i$ such that 
\begin{equation*}
  \rho_*([\widetilde{M}]) = b_0 + b_1(1-t) + \cdots + b_{n-1}(1-t)^{n-1}.
\end{equation*}
To prove the lemma, it suffices to show
\begin{equation} \label{E12}
\begin{aligned} 
& b_0  = \cdots = b_{n-m-1} = 0   \quad \text{ and} \\
& a_i = b_{n-m+i}  \quad \text{  for each $i = 0, \dots, m-1$}.
\end{aligned}
\end{equation}

Recall $q \colon X~\to \Spec{\sk}$ and $s \colon~\bP^{n-1} \to
\Spec{\sk}$ are the structure maps; see (\ref{spaces}).  The class of a
coherent sheaf $\cF$ on $X$ maps to $\sum_i (-1)^i \dm_k H^i(X, \cF)$
under $q_*$, and similarly for $s_*$.  The sheaf cohomology of the 
line bundles $\cO_{\bP^{n-1}}(i)$ \cite[Exercise III.5.4]{Ha} shows that $s_*$
sends each $(1-t)^i \in \bQ[t]/(1-t)^n \iso \K{\bP^{n-1}}_\bQ$ to $1$, 
for $i = 0, \dots, n-1$, and so
\begin{equation*}
  s_* \big(\rho_*([\widetilde{M}])(1-t)^i\big) = b_0 + \cdots + b_{n-i-1}.
\end{equation*}

On the other hand, for a finitely generated graded $R$-module $T$, we
have $q_*([\widetilde{T}]) = \Hp[T]{0}$ and
hence $q_*([\widetilde{T}(i)]) = \Hp[T]{i}$ for all $i$ \cite[Exercise III.5.2]{Ha}.  It follows
that
\begin{equation*}
q_*([\widetilde{T}](1 - [\cO(-1)])) = 
\Hp[T]{0} - \Hp[T]{-1} = 
\Hpi[T]{1}{0}.
\end{equation*}
 From this we deduce that
$$
s_* \big(\rho_*([\widetilde{M}])(1-t)^i\big) = q_*\big([\widetilde{M}](1 - [\cO(-1)])^i\big)
= \Hpi[M]{i}{0}.
$$
Using \eqref{E11}, we conclude that 
$$
a_0 + \cdots + a_{m -i - 1} = b_0 + \cdots b_{n-i-1}
$$
for all $i \geq 0$.
The equations \eqref{E12} follow.
\end{proof}

The preceding lemma relates $\rho_*([\widetilde{M}])$ to an invariant that
is closely related to the Hilbert polynomial of $M$. The next one
relates $\rho_*$ to $\theta$.

\begin{lem} \label{rho and theta} Let $X$ be as in \eqref{assumptions} with $\sk$ infinite.  Then for any pair of elements $x,y \in \K{X}_\bQ$, there is an equality
\begin{equation*}
\frac{(1-t)^{n-1}}{d^2} \theta(x,y)
= 
\left( \frac{\rho_*(x)}{\rho_*(1)} \right)
\left( \frac{\rho_*(y)}{\rho_*(1)} \right) 
- \left(\frac{\rho_*(x\cdot y)}{\rho_*(1)} \right) 
\end{equation*}
in $\K{\bP^{n-1}}_\bQ \iso \bQ[t]/(1-t)^n$.

Consequently,
\begin{equation*}
\theta(x,y) =  s_*\left[
\left(\frac{d \cdot \rho_*(x)}{\rho_*(1)} \right)
\left(\frac{d \cdot \rho_*(y)}{\rho_*(1)} \right) 
- d\left(\frac{d \cdot \rho_*(x\cdot y)}{\rho_*(1)} \right) \right] 
\end{equation*}
where $s: \bP^{n-1} \to \Spec{\sk}$ is the structure map. 
\end{lem}

\begin{proof}
  The second equation follows from the first since $s_*((1-t)^{n-1}) =
  1$.  In this proof, we suppress the variable from the notation
  $\Hs[M]{t}$ and simply write
  \renewcommand{\Hs}[2][]{\operatorname{H}_{#1}}$\Hs[M]{t}$.

  As $\theta$ is bilinear and $\rho_*$ is linear, we may assume $x =
  [\widetilde{M}]$ and $y = [\widetilde{N}]$ for finitely generated
  graded $R$-modules $M$ and $N$. Let $\Hs[i]{t} =
  \Hs[\Tor_i^R(M,N)]{t}$, the Hilbert series of $\Tor_i^R(M,N)$. By
  \cite[Lemma 7]{AB}, we have
\begin{equation} \label{E1}
\sum_{i \geq 0} (-1)^i \Hs[i]{t} = \frac{\Hs[M]{t} \Hs[N]{t}}{\Hs[R]{t}}.
\end{equation}

For a sufficiently large even integer $E$, the
length of $\Tor_i^R(M,N)$ is finite and there is an isomorphism of
graded $R$-modules
\begin{equation*}
  \Tor_i^R(M,N)(-d) \iso \Tor_{i+2}^R(M,N), \qquad {\text{for each }} i \geq E.
\end{equation*}
As these torsion modules are finite length, and hence are nonzero in
only finitely many graded degrees, it follows that $\Hs[E]{t}$ and
$\Hs[E+1]{t}$ are polynomials.  Moreover, 
\begin{equation*}
\Hs[E+2j]{t} = t^{dj} \Hs[E]{t} 
\quad \text{and} \quad 
\Hs[E+1+2j]{t} = t^{dj} \Hs[E+1]{t},
\quad \text{for all $j \geq 0$}.
\end{equation*}
Consequently 
\begin{equation*}
\begin{split}
\sum_{i \geq E} (-1)^i \Hs[i]{t} 
 = (\Hs[E]{t} - \Hs[E+1]{t})(1 + t^d + t^{2d} + \cdots)
 = \frac{\Hs[E]{t} - \Hs[E+1]{t}}{1-t^d}
= \frac{\Hs[E]{t} - \Hs[E+1]{t}}{e_R(t)(1-t)}.
\end{split}
\end{equation*}
Combining this with \eqref{E1} gives
\begin{equation*}
\sum_{i=0}^{E-1} (-1)^i \Hs[i]{t} + 
\frac{\Hs[E]{t} - \Hs[E+1]{t}}{e_R(t)(1-t)}
=
\frac{\Hs[M]{t} \Hs[N]{t}}{\Hs[R]{t}}
= 
\frac{(1-t)^n\Hs[M]{t} \Hs[N]{t}}{e_R(t)}.
\end{equation*}
Multiplying both sides by $(1-t)^n/e_R(t)$ and
rearranging the terms gives
\begin{equation} \label{E2}
\frac{(1-t)^n\Hs[M]{t}}{e_R(t)} 
\frac{(1-t)^n\Hs[N]{t}}{e_R(t)} -
\sum_{i=0}^{E-1} (-1)^i \frac{(1-t)^n}{e_R(t)} \Hs[i]{t} =
\frac{(\Hs[E]{t} - \Hs[E+1]{t})}{(e_R(t))^2}(1-t)^{n-1}.
\end{equation}
Both sides of this equation are power series in powers of $1-t$ and we
may thus take their images in $\bQ[t]/(1-t)^n$. We claim doing so
results in the equation in the statement of this lemma.  For a finitely generated graded $R$-module $T$, \lemref{rho and HS}
shows that 
$$ 
\rho_*([\widetilde{T}])/\rho_*(1) =
(1-t)^n\Hs[T]{t}/e_R(t)
\quad \text{in $\bQ[t]/(1-t)^n$.}
$$
Since the coherent
$\cO_X$-sheaf associated to a graded module of finite length is zero,
we have $\widetilde{\Tor_i^R(M,N)} = 0 $ for all $i \geq E$.  Hence,
in the ring $\K{X}_\bQ$ we have
\begin{equation*}
x \cdot y = \sum_{i=0}^{E-1} (-1)^i [\widetilde{\Tor_i^R(M,N)}].  
\end{equation*}
The image of the left-hand side of \eqref{E2} in the ring
$\bQ[t]/(1-t)^n$ is therefore
\begin{equation*}
\frac{\rho_*(x)}{\rho_*(1)} \frac{\rho_*(y)}{\rho_*(1)} 
- \sum_{i=0}^{E-1} (-1)^i \frac{\rho_*([\widetilde{\Tor_i^R(M,N)}])}{\rho_*(1)}
=
\frac{\rho_*(x)}{\rho_*(1)}\frac{\rho_*(y)}{\rho_*(1)} - 
\frac{\rho_*(x \cdot y)}{\rho_*(1)}.
\end{equation*}
To simplify the right-hand side of \eqref{E2}, observe that $f(t) =
(\Hs[E]{t} - \Hs[E+1]{t})/(e_R(t))^2$ is a rational function without a
pole at $t = 1$.  Modulo $(1 - t)^n$, we have 
\begin{equation*}
 f(t)(1 - t)^{n-1} = f(1)(1-t)^{n-1} + \frac{f(t) - f(1)}{1 -
   t}(1-t)^n \equiv f(1)(1-t)^{n-1}.
\end{equation*}
Since $\Hs[E][(1)= \length({\Tor_E^R(M,N)})$, $\Hs[E+1][(1) = \length({\Tor_{E+1}^R(M,N)})$, and $e_R(1) = d$, it follows that
\begin{equation*}
\frac{\rho_*(x)}{\rho_*(1)}\frac{\rho_*(y)}{\rho_*(1)} - 
\frac{\rho_*(x \cdot y)}{\rho_*(1)} = f(1)(1 - t)^{n-1} = 
\frac{\theta(x,y)}{d^2}(1 - t)^{n-1}.  
\qedhere
\end{equation*}

\renewcommand{\Hs}[2][]{\operatorname{H}_{#1}(#2)}
\end{proof}

\begin{lem} \label{rr}
Let $X$ be as in \eqref{assumptions} with $k$ an infinite field.  The diagram
\begin{equation} \label{ECH}
\begin{gathered}
\xymatrix{
\K{X}_\bQ \ar[rr]^-{\ds \frac{d}{\rho_*(1)} \cdot \rho_*} 
         \ar[d]_-{\ds \ch{}} 
  & & \K{\bP^{n-1}}_\bQ \ar[d]^-{\ds \ch{}} \\
CH^\bu(X)_\bQ \ar[rr]_-{\ds\rho_*} 
  & & CH^\bu(\bP^{n-1})_\bQ
}
\end{gathered}
\end{equation}
commutes.

If $k$ is separably closed, then 
$\rho_* \circ \ch[\et] =  {\ds \frac{d}{\rho_*(1)}} \cdot  \ch[\et] \circ \rho_*$. 

If $\sk = \C$, then $\rho_* \circ \ch[\topo] = {\ds \frac{d}{\rho_*(1)}} \cdot
\ch[\topo] \circ \rho_*$.
\end{lem}

\begin{proof} 
  The  maps $cy_{\et}$ and $cy_{\topo}$ combine with
  \eqref{ECH} to give the last results since they commute with push-forwards.

By the Grothendieck-Riemann-Roch Theorem \cite[15.2]{Fu},
for any $x \in \K{X}_\bQ$, 
\begin{equation*}
\rho_*(\ch{x} \cup \Td{\cT_X}) =
\ch{\rho_*(x)} \cup \Td{\cT_{\bP^{n-1}}}
\end{equation*}
where $\cT_Y$ is the tangent bundle of a smooth variety $Y$.  Now
\begin{equation*}
\ch{\left(\frac{d \cdot \rho_*(x)}{\rho_*(1)}\right)}
= \frac{d \cdot \ch{\rho_*(x)}}{\ch{\rho_*(1)}}
= \frac{d \cdot \ch{\rho_*(x)} \cup \Td{\cT_{\bP^{n-1}}}}
               {\ch{\rho_*(1)} \cup \Td{\cT_{\bP^{n-1}}}}
= \frac{d \cdot \rho_*\left(\ch{x} \cup \Td{\cT_X}\right)}
{\rho_*\left(\ch{1} \cup \Td{\cT_X}\right)}.
\end{equation*}
But $\Td{\cT_X} \in 1 + \beta CH^\bu(X)_\bQ$ as per Lemma \ref{ch and tau agree} and hence $\Td{\cT_X} =
\rho^*(\xi)$ for some unit $\xi \in \CH[\bullet]{\bP^{n-1}}_\bQ$.
The projection formula gives
\begin{equation*}
\frac{d \cdot \rho_*\left(\ch{x} \cup \Td{\cT_X}\right)}
{\rho_*\left(\ch{1} \cup \Td{\cT_X}\right)}=
\frac{d \cdot \rho_*\left(\ch{x} \cup \rho^*(\xi)\right)}
{\rho_*\left(\ch{1} \cup \rho^*(\xi)\right)} =
\frac{d \cdot \rho_*(\ch{x}) \xi}{\rho_*(\ch{1}) \xi}.
\end{equation*}
Since
$\rho_*(\ch{1}) = d \in CH^0(\bP^{n-1})_\bQ$, the diagram \eqref{ECH} commutes.
\end{proof}

\section{Proofs of the main theorems} \label{proof section}

This section contains the proofs of the main results of this paper,
which are stated in Section \ref{statement section}. We use the
results developed in Section \ref{Hilbert section}.

\begin{proof}[Proof of Proposition \ref{thm}]
  Multiplying the first equation in \lemref{rho and theta} by $d^2$, applying $\ch[\et]$, and then simplifying using
  \lemref{rr} yields
\begin{equation*}
  \rho_*(\ch[\et]{x})\cup\rho_*(\ch[\et]{y})
  -
  d \cdot \rho_*(\ch[\et]{(x\cdot y)})
  =
  \theta(x,y) \ch[\et]{(1-t)^{n-1}}
\end{equation*}
for all $x,y \in \K{X}_\bQ$.  As $t$ corresponds to
$[\cO_{\bP^{n-1}}(-1)]$ under the isomorphism
$\K{\bP^{n-1}}_\bQ \iso \bQ[t]/(1 - t)^n$, the \'etale Chern
character of $t$ is $\ch[\et]{t} = \exp(-\varsigma) = 1 - \varsigma +
\frac{\varsigma^2}{2} - \frac{\varsigma^3}{3!}  + \cdots$ \cite[Example 3.2.3]{Fu}, where
$\varsigma \in \H[\et]{2}{\bP^{n-1}}{\bQl(1)}$ is the class of a
hyperplane.  Since $\varsigma^n = 0$ in
$\H[\et]{ev}{\bP^{n-1}}{\bQl}$, it follows that $\ch[\et](1 -
t)^{n-1}$ is $\varsigma^{n-1}$.
Using this gives
$$ 
  \int_{\bP^{n-1}} \Big(
  \rho_*(\ch[\et]{x})\cup\rho_*(\ch[\et]{y}) - d \cdot
  \rho_*\big(\ch[\et]{(x\cdot y)}\big) \Big) 
  = \theta(x,y) \int_{\bP^{n-1}} \varsigma^{n-1}.
$$
The first
assertion of Proposition \ref{thm} follows from the fact that
$\int_{\bP^{n-1}} \varsigma^{n-1} = 1$.

For a proof of the second assertion, use $ch_\topo$ in place of
$ch_\et$.
\end{proof}

\begin{proof}[Proof of Theorem \ref{n even result}] %
  We may assume $\sk$ is separably closed by applying \eqref{theta
    field stable} with $k' = \sk_{sep}$.  Define a pairing $\phi$ on
  $H_{\et}^{ev}(X, \bQl)$ by the formula
\begin{equation} \label{E821b}
\phi(a,b) = \int_{\bP^{n-1}} \big( \rho_*(a)\cup\rho_*(b) - d \cdot
  \rho_*(a \cup b) \big).
\end{equation}
Since $\ch[\et]$ is a ring homomorphism, Proposition \ref{thm} shows that for all
$x, y$  in $K(X)_\bQ$, $\theta(x,y) = \phi(\ch[\et] x, \ch[\et] y)$.
Thus \thmref{n even result} follows
from the assertion that $\phi = 0$ when $n$ is even.

Using the projection formula, we get that $\phi(a,b) = 0$ if either
$a$ or $b$ lies in the image of the ring map $\rho^* \colon
\H[\et]{ev}{\bP^{n-1}}{\bQl} \to \H[\et]{ev}{X}{\bQl}$. Recall that $\gamma = \cy[\et](\beta)
\in H^2_\et(X, \bQl(1))$ is the \'etale cohomology
class of the divisor given by a generic hyperplane section of $X$.
Since $\gamma$
lies in the image of $\rho^*$, we have $\phi(a,b) = 0$ if either $a$ or
$b$ is a multiple of a power of $\gamma$. The theorem therefore
follows from \eqref{etale cohomology}.
\end{proof}

\begin{proof}[Proof of \thmref{n odd thm}]
 We use the pairing $\phi$ introduced in
 \eqref{E821b}. Using \eqref{etale cohomology} and the fact that $\phi(a,b)
 = 0$ if either $a$ or $b$ is a multiple of a power of $\gamma$, we
 get that the pairing $\phi$ factors through the
 canonical surjection
\begin{equation*}
\H[\et]{ev}{X}{\bQl} \onto
\frac{\H[\et]{n-1}{X}{\bQl(\frac{n-1}{2})}}{\bQl \cdot \gamma^{\frac{n-1}{2}}}
\end{equation*}
when $n$ is odd.  We claim
\begin{equation} \label{E914}
\sint_{\bP^{n-1}} (\rho_*(a)\cup\rho_*(b))
=
\sint_X (a \cup \gamma^{\frac{n-1}{2}})
\sint_X (b \cup \gamma^{\frac{n-1}{2}}),
\end{equation}
for all $a, b \in H^{n-1}_{\et}(X, \bQl(\frac{n-1}{2}))$.
To see this, first recall that $H^{ev}_{\et}(\bP^{n-1}, \bQl) \cong
\bQl[\varsigma]/\ip{\varsigma^n}$ and that
$$
\sint_{\bP^{n-1}}
\varsigma^i =  \begin{cases}
0, & \text{for $i = 0, \cdots, n-2$ and} \\
1, & \text{for $i = n-1$.} \\
\end{cases}
$$
If $\rho_*(a) = r \varsigma^{\frac{n-1}{2}}$
and $\rho_*(b) = r' \varsigma^{\frac{n-1}{2}}$ for $r,r' \in \bQl$,
then $\sint_{\bP^{n-1}} (\rho_*(a) \cup \rho_*(b)) = rr'$.
On the other hand,
since $\rho$ factors as $X \into \bP^n \setminus
\{Q\} \onto \bP^{n-1}$ (with the second map being linear projection
away from $Q$), we have $\rho^*(\varsigma) = \gamma$. The equation
$\sint_X = \sint_{\bP^{n-1}} \circ \rho_*$ and
the
projection formula give
$$
\sint_X(a \cup \gamma^{\frac{n-1}{2}}) =
\sint_{\bP^{n-1}} \rho_*(a \cup \gamma^{\frac{n-1}{2}}) =
\sint_{\bP^{n-1}} \rho_*(a)  \cup \varsigma^{\frac{n-1}{2}} = r.
$$
Similarly, $\sint_X(b \cup \gamma^{\frac{n-1}{2}}) = r'$, and
\eqref{E914} follows.

Since \eqref{E914} holds, formula \eqref{E821b} becomes
$$
\phi(a,b) =
\sint_X (a \cup \gamma^{\frac{n-1}{2}})
\sint_X (b \cup \gamma^{\frac{n-1}{2}})
- d \cdot \sint_X (a \cup b).
$$
The first assertion of the theorem now follows from Proposition \ref{thm}.

The proof of the
second assertion is analogous.
\end{proof}

\begin{proof}[Proof of Theorem \ref{n odd pos defn}] %
  Recall that for integers $p,q$ with $p + q = n - 1$, the complex
  vector space $H^{p,q}(X(\bC))$ is the $(p,q)$-part of the Hodge
  decomposition of $\H{n-1}{X(\bC)}{\bQ}$.  
(See \cite[\S 0.6]{GH} or \cite[19.3.6]{Fu}.)
  Define 
\begin{equation*}
W =   \H{n-1}{X(\bC)}{\bQ} \, \cap \, H^{\frac{n-1}{2},\frac{n-1}{2}}(X(\bC)),
\end{equation*}
a $\bQ$-vector subspace of $\H{n-1}{X(\bC)}{\bQ}$.
Since the image of $cy_{\topo}^{\frac{n-1}{2}}$ is contained in $W$ \cite[19.3.6]{Fu},
it follows that the image of $\ch[\topo]^{\frac{n-1}{2}} \colon \K{X}_\bQ \to 
\H{n-1}{X(\bC)}{\bQ}$ is also contained in $W$.
We argue that the restriction of 
$(-1)^{\frac{n+1}{2}}\theta_\topo$ to $W/\bQ \cdot \gamma^{\frac{n-1}{2}}$
is positive definite.   

Define an injection 
$e \colon W / \bQ \cdot \gamma^{\frac{n-1}{2}} \into
\H{n-1}{X(\bC)}{\bQ}$ by setting, for $a \in W$,
\begin{equation*}
  e(a) = %
  a - \frac{\sint_X a \cup \gamma^{\frac{n-1}{2}} }{d} \gamma^{\frac{n-1}{2}}
  \in \H{n-1}{X(\bC)}{\bQ}.  
\end{equation*}
The image of $e$ is contained in $H^{\frac{n-1}{2},
  \frac{n-1}{2}}(X(\C))$ since both $W$ and $\gamma^{\frac{n-1}{2}}$
are.  It is also contained in the primitive part of $H^{n-1}(X(\bC),
\bQ)$, i.e., in the subspace of elements $v$ with $\gamma \cup v = 0$.
Indeed, \eqref{etale cohomology} implies that repeated
cupping with $\gamma$, followed by the map $\sint_X$, forms a sequence of
isomorphisms:
\begin{equation*}
  \xymatrix{
    \H{n+1}{X(\C)}{\bQ} 
    \ar[r]^-{\gamma \cup - }_-\iso 
    & \H{n+3}{X(\C)}{\bQ} 
      \ar[r]^-{\gamma \cup - }_-\iso  
      & \cdots 
        \ar[r]^-{\gamma \cup - }_-\iso
        & \H{2n-2}{X(\C)}{\bQ} 
          \ar[r]^-{\sint_X}_-\iso 
          & \bQ.
  }
\end{equation*}
So the vanishing of $\gamma \cup e(a)$ follows from the vanishing of
$\int_X \gamma^{\frac{n-1}{2}} \cup e(a)$, which is clear.

Let $Q$ be the bilinear pairing on $\H{n-1}{X(\bC)}{\bQ}$ given by 
$
  Q(x,y) = \sint_X x \cup y.  
$
Straightforward computation verifies that for $a, b \in W$, 
\begin{equation*}
 \theta_\topo(a,b) = - d \cdot Q(e(a), e(b)).
\end{equation*}
The Hodge-Riemann bilinear relations \cite[page 165]{Vo} give that
$(-1)^{\frac{(n-1)(n-2)}{2}} Q$ is positive definite on the primitive
part of 
$W$ and hence on the
image of $e$.  

For the last assertion of the theorem, let $R$ be as in
\eqref{assumptions} with $\sk$ a field of characteristic zero, and
let $M$ be a finitely generated $R$-module.  There is a finitely
generated field extension of $\bQ$ that contains the coefficients
of $f(x_0, \dots, x_n)$ and the entries of a presentation matrix for $M$.
Using \eqref{theta field stable} twice, we may assume $\sk \subseteq
\bC$ and then $\sk = \bC$.  The result follows from the first
assertion of the theorem. 
\end{proof}

\section{Generalization of the Main Theorems}  \label{sec:nonstandardgrading}

The following theorem allows us, among other things, to extend our
result to hypersurfaces that are homogeneous with respect to a non-standard grading.

\begin{thm} \label{nonhom} Let $R$ be as in
  \eqref{assumptions}. Suppose $S \subset R$ is a subring such that
  the the map $S \into R$ is finite, flat, and a local complete
  intersection (see \cite[B.7.6]{Fu}).
Note that $n = \dim R = \dim S$.  Assume
  $S$ is also a hypersurface with isolated singularity.
\begin{enumerate}
\item %
  If $n$ is even and $\sk$ is an arbitrary field, then $\theta^S$
  vanishes on $\G{S}_\bQ$.
\item 
  If $n$ is odd and $\sk$ is an arbitrary field, then $\theta^S$ on
  $\G{S}_\bQ$ is induced from $\theta_{CH(S)}^S$ on
  $\CH[\frac{n-1}{2}]{S}_\bQ$ via the surjective map
  $\tau^{\frac{n-1}{2}} \!\colon \!\G{S}_\bQ \onto
  \CH[\frac{n-1}{2}]{S}_\bQ$.
\item %
  If $n$ is odd and $\chr(k) = 0$, then $(-1)^{\frac{n+1}{2}}\theta^S$ is positive semi-definite
 on $\G{S}_\bQ$. 
 \item %
  If $n = 3$ and $\sk = \bC$, then $\theta_{CH(S)}^S$ is positive definite on
  $\CH[1]{S}_\bQ$. In particular, $\theta^S = 0$
  if and only if $\CH[1]{S}$ is torsion.
\end{enumerate}
\end{thm}

\begin{proof}
  Let $h \colon \Spec{R} \to \Spec{S}$ be the induced map.  Since $h_*
  \circ h^*$ is multiplication by the degree of $h$, the map
 $ h^* $ 
 is injective.  By \cite[18.2]{Fu} the diagram
\begin{equation*} 
  \xymatrix{
    \G{S}_\bQ 
    \ar[r]^-{\ds \tau}_-{\ds \iso} 
    \ar[d]_-{\ds h^*} 
    & \CH[\bullet]{S}_\bQ
    \ar[d]^-{\ds \Td{T_h} h^*} 
    \\
    \G{R} 
    \ar[r]_-{\ds \tau}^-{\ds \iso}
    & \CH[\bullet]{R}_\bQ
  }
\end{equation*} 
commutes, where $T_h$ is the virtual tangent bundle of $h$.  Let
$\mathfrak n$ denote the isolated singularity of $S$ and let $c =
\length_R(R/\mathfrak n R)$.  Then the following triangle commutes.
\begin{equation*}
  \xymatrix@=1.25em{
    \G{S}_\bQ^{\otimes 2} 
    \ar[dr]_(.5){\ds \theta^S} 
    \ar[rr] 
    & & \G{R}_\bQ^{\otimes 2} \ar[ld]^(.5){\ds \frac{1}{c} \theta^{R}} \\
    & \bQ
  }
\end{equation*}

The first assertion follows from this commutative triangle and Theorem
\ref{n even result}.

The second assertion follows from the above two commutative diagrams
and Corollary \ref{theta on Chow R} by taking
\begin{equation*}
\theta_{CH}^S = \frac{1}{c} \cdot \theta_{CH}^R \circ (\Td{T_h}
h^*)^{\tensor 2}.
\end{equation*}

The last two assertions follow from Theorem \ref{n odd pos defn} and
Corollary \ref{n odd theta = 0}, using the fact that $\Td{T_h} h^*$
is injective since $\Td{T_h}$ is a unit.
\end{proof}

We now describe a general situation to which Theorem \ref{nonhom}
applies.

\begin{example} \label{German hypersurface} %
  Let $S= \sk[y_0, \dots, y_n]/(g)$ be graded with $\deg(y_i) = e_i$
  for some integers $e_i \geq 1$, and let $g(y_0, \dots, y_n)$ be
  homogeneous of degree $d \geq 1$.  Set $R = \sk[x_0, \dots,
  x_n]/(f)$, where $\deg(x_i) = 1$ and $f = g(x_0^{e_0}, \dots,
  x_n^{e_n})$. Then $S$ is a subring of $R$, with $y_i = x_i^{e_i}$.

  Since $R \iso S[x_0, \dots, x_n]/\ip{x_i^{e_i} - y_i\, | \, i = 0,
    \ldots, n}$, it follows that the map $S \into R$ is finite, flat,
  and a local complete intersection.  When both $S$ and $R$ have
  isolated singularities, Theorem \ref{nonhom} applies; in
  particular, it applies to rings of the form
\begin{equation} \label{general S}
S = \sk[y_0, \dots, y_n]/(y_0^{m_0} + \cdots +
y_n^{m_n}), \qquad {\text{with $m_i \geq 1$ for all $i$.}}
\end{equation}
\end{example}

\begin{example} (A specific application of Example \ref{German hypersurface}.)
J.~Bingener and U.~Storch \cite[\S 12]{BS} show that for $S$ given by
\eqref{general S} with $n=3$ and $\sk$ algebraically closed of characteristic 
zero, the group $\CH[1]{S}$ is finitely
generated free abelian.  They find bounds on its rank.  
For instance,
if $S$ is defined by the 4-tuple $(m_0, m_1, m_2, m_3) = (2, 3, 3,
6)$, then $\CH[1]{S} \iso \bZ^4$, and therefore $\theta^S \neq 0$.  
In this case, the module $M = S/(y_0-y_3^3, y_1+y_2)$ determines a
non-zero class in $CH^1(S)_\Q$  and hence
$\theta^S(M,M) \neq 0$.  On the other hand, if $m_3 = 7$ instead, then $\CH[1]{S} =0$, and
hence $\theta^S$ vanishes.

When $m_0, m_1,
m_2$ are distinct primes and $m_3 = m_0m_1m_2$,
Bingener and Storch find an upper bound on the
rank of $\CH[1]{S}$; see
\cite[(12.9)(2)]{BS}.  For most primes, the exact rank of $\CH[1]{S}$
is unknown, and $\CH[1]{S}$ may even be trivial. For
$$
S = \sk[y_0, \dots, y_3]/(y_0^2 + y_1^3+ y_2^5 +
y_3^{30}),
$$ 
they show that $\CH[1]{S} \iso \bZ^8$, and therefore
$\theta^S \neq 0$. It appears difficult to find generators of
$\CH[1]{S}$. Indeed, we were unable to find a single explicit module
$M$ for which $\theta^S(M,M) \ne 0$, even though our results show such
modules must exist.
\end{example}

\bibliographystyle{amsplain}

\begin{thebibliography}{10}

\bibitem{SGA4V1}
\emph{Th\'eorie des topos et cohomologie \'etale des sch\'emas. {T}ome 1:
  {T}h\'eorie des topos}, Lecture Notes in Mathematics, Vol. 269,
  Springer-Verlag, Berlin, 1972, S{\'e}minaire de G{\'e}om{\'e}trie
  Alg{\'e}brique du Bois-Marie 1963--1964 (SGA 4), Dirig{\'e} par M. Artin, A.
  Grothendieck, et J. L. Verdier. Avec la collaboration de N. Bourbaki, P.
  Deligne et B. Saint-Donat.

\bibitem{SGA4V2}
\emph{Th\'eorie des topos et cohomologie \'etale des sch\'emas. {T}ome 2},
  Lecture Notes in Mathematics, Vol. 270, Springer-Verlag, Berlin, 1972,
  S{\'e}minaire de G{\'e}om{\'e}trie Alg{\'e}brique du Bois-Marie 1963--1964
  (SGA 4), Dirig{\'e} par M. Artin, A. Grothendieck et J. L. Verdier. Avec la
  collaboration de N. Bourbaki, P. Deligne et B. Saint-Donat.

\bibitem{SGA4V3}
\emph{Th\'eorie des topos et cohomologie \'etale des sch\'emas. {T}ome 3},
  Lecture Notes in Mathematics, Vol. 305, Springer-Verlag, Berlin, 1973,
  S{\'e}minaire de G{\'e}om{\'e}trie Alg{\'e}brique du Bois-Marie 1963--1964
  (SGA 4), Dirig{\'e} par M. Artin, A. Grothendieck et J. L. Verdier. Avec la
  collaboration de P. Deligne et B. Saint-Donat.

\bibitem{SGA5}
\emph{Cohomologie {$l$}-adique et fonctions {$L$}}, Lecture Notes in
  Mathematics, Vol. 589, Springer-Verlag, Berlin, 1977, S{\'e}minaire de
  G{\'e}ometrie Alg{\'e}brique du Bois-Marie 1965--1966 (SGA 5), Edit{\'e} par
  Luc Illusie.

\bibitem{AB}
Luchezar~L. Avramov and Ragnar-Olaf Buchweitz, \emph{Lower bounds for {B}etti
  numbers}, Compositio Math. \textbf{86} (1993), no.~2, 147--158.

\bibitem{BS}
J{\"u}rgen Bingener and Uwe Storch, \emph{Zur {B}erechnung der
  {D}ivisorenklassengruppen kompletter lokaler {R}inge}, Nova Acta Leopoldina
  (N.F.) \textbf{52} (1981), no.~240, 7--63, Leopoldina Symposium:
  Singularities (Th{\"u}ringen, 1978).

\bibitem{Bu}
Ragnar-Olaf Buchweitz, \emph{Maximal {C}ohen-{M}acaulay modules and {T}ate
  cohomology over {G}orenstein rings}, Preprint, 1966.

\bibitem{CE}
Henri Cartan and Samuel Eilenberg, \emph{Homological algebra}, Princeton
  University Press, Princeton, New Jersey, 1956.

\bibitem{CHWW}
Guillermo Corti{\~n}as, Christian Haesemeyer, Mark~E. Walker, and Chuck Weibel,
  \emph{The {$K$}-theory of toric varieties}, Trans. Amer. Math. Soc.
  \textbf{361} (2009), 3325--3341.

\bibitem{Dao1}
Hailong Dao, \emph{Decency and rigidity over hypersurfaces},
  arXiv:math/0611568, 2008.

\bibitem{Dao2}
\bysame, \emph{Some observations on local and projective hypersurfaces}, Math.
  Res. Lett. \textbf{15} (2008), no.~2, 207--219.

\bibitem{SGA4.5}
P.~Deligne, \emph{Cohomologie \'etale}, Lecture Notes in Mathematics, Vol. 569,
  Springer-Verlag, Berlin, 1977, S{\'e}minaire de G{\'e}om{\'e}trie
  Alg{\'e}brique du Bois-Marie SGA 4$\frac{1}{2}$, Avec la collaboration de J.
  F. Boutot, A. Grothendieck, L. Illusie et J. L. Verdier.

\bibitem{FK}
Eberhard Freitag and Reinhardt Kiehl, \emph{\'{E}tale cohomology and the {W}eil
  conjecture}, Ergebnisse der Mathematik und ihrer Grenzgebiete (3) [Results in
  Mathematics and Related Areas (3)], vol.~13, Springer-Verlag, Berlin, 1988,
  Translated from the German by Betty S. Waterhouse and William C. Waterhouse,
  With an historical introduction by J. A. Dieudonn{\'e}.

\bibitem{Fu}
William Fulton, \emph{Intersection theory}, Ergebnisse der Mathematik und ihrer
  Grenzgebiete (3) [Results in Mathematics and Related Areas (3)], vol.~2,
  Springer-Verlag, Berlin, 1984.

\bibitem{GH}
Phillip Griffiths and Joseph Harris, \emph{Principles of algebraic geometry},
  Wiley-Interscience [John Wiley \& Sons], New York, 1978, Pure and Applied
  Mathematics.

\bibitem{Ha}
Robin Hartshorne, \emph{Algebraic geometry}, Springer-Verlag, New York, 1977,
  Graduate Texts in Mathematics, No. 52.

\bibitem{Heit}
Raymond Heitmann, \emph{The direct summand conjecture in dimension three}, Ann.
  of Math. \textbf{156} (2002), no.~2, 695--712.

\bibitem{Ho2}
Melvin Hochster, \emph{Contracted ideals from integral extensions of regular
  rings}, Nagoya Math. J. \textbf{51} (1973), 25--43.

\bibitem{Ho}
\bysame, \emph{The dimension of an intersection in an ambient hypersurface},
  Algebraic geometry ({C}hicago, {I}ll., 1980), Lecture Notes in Math., vol.
  862, Springer, Berlin, 1981, pp.~93--106.

\bibitem{Kl}
Steven~L. Kleiman, \emph{The standard conjectures}, Motives ({S}eattle, {WA},
  1991), Proc. Sympos. Pure Math., vol.~55, Amer. Math. Soc., Providence, RI,
  1994, pp.~3--20.

\bibitem{Se}
Jean-Pierre Serre, \emph{Alg\`ebre locale. {M}ultiplicit\'es}, Cours au
  Coll\`ege de France, 1957--1958, r\'edig\'e par Pierre Gabriel. Seconde
  \'edition, 1965. Lecture Notes in Mathematics, vol.~11, Springer-Verlag,
  Berlin, 1965.

\bibitem{Vo}
Claire Voisin, \emph{Transcendental methods in the study of algebraic cycles},
  Algebraic cycles and {H}odge theory ({T}orino, 1993), Lecture Notes in Math.,
  vol. 1594, Springer, Berlin, 1994, pp.~153--222.

\end{thebibliography}

\providecommand{\bysame}{\leavevmode\hbox to3em{\hrulefill}\thinspace}
\providecommand{\MR}{\relax\ifhmode\unskip\space\fi MR }
\providecommand{\MRhref}[2]{%
  \href{http://www.ams.org/mathscinet-getitem?mr=#1}{#2}
}
\providecommand{\href}[2]{#2}

\end{document}